\documentclass[11pt]{amsart}
\usepackage{hyperref}

\begin{document}
\title{Beyond the Black Box}
\author{Jeroen Demeyer, William Stein, and Ursula Whitcher}
\maketitle

``Commercial computer algebra systems are black boxes, and their algorithms are opaque to the users,''  complained a trio of mathematicians whose ``misfortunes'' are detailed in a recent \emph{Notices} article, \cite{trio}.  ``We reported the bug on October 7, 2013 \dots By June 2014, nothing had changed \dots All we could do was wait.''

In open source software, the code underlying each operation is available to anyone who chooses to look at it.  Because open-source code can be checked directly, it is highly valuable for replicable, peer-reviewable research.  We are users and developers of SageMath, open-source software by and for the mathematical community (see \cite{william}).  SageMath is a full computer algebra system that can be installed locally or accessed freely in a web browser through the SageMathCloud, at \cite{cloud}.  Sage also provides a consistent interface to many free libraries of mathematical software.

Open source is not a panacea: errors can arise in open-source software as well as closed-source software.  We trace the history of a representative bug in Sage, to illustrate the role the mathematical community plays in detecting and fixing bugs in open-source software.  The unfortunate trio \cite{trio} found a bug in \emph{Mathematica}'s algorithm for determinants of large integer matrices.  SageMath evaluates that determinant correctly (see \cite{worksheet}).  For the sake of comparison, we describe a different past problem with Sage's computation of determinants, and the process followed in resolving it.

The first step in fixing a bug is that somebody must notice a problem.  In our case, SageMath release 5.6 stalled when computing the determinant of a large integer matrix: the computation never completed.

Next, the person who notices the bug must report it.  Typically, a user would post a question to the sage-support or sage-devel Google group.  Asking questions simultaneously warns active SageMath users about the problem and advertises for developers who might fix it.  Bugs can also be reported by emailing \texttt{bugs@sagemath.com}.  A very confident user might move straight to the next step: creating a ticket.

Proposed changes to the SageMath code are tracked on the SageMath trac server, \cite{trac}.  Individual bugs or code enhancements are called tickets.  In January 2013, our determinant bug became ticket 14032, with the initial title \emph{determinant() of large integer matrices broken} (see \cite{14032}).  The person who creates a ticket chooses a priority for the fix, from the options of ``trivial'', ``minor'', ``major'', ``critical'', and ``blocker''.  Most SageMath tickets are given the middle ranking of ``major''; our determinant bug was ``critical''.  A new version of SageMath cannot be released if there are open tickets marked ``blocker''.

Once the ticket is created, developers can start hunting down the bug's origins.  In our case, the first clue was that stalling arose when computing determinants of integer matrices of size greater than $50 \times 50$.  This was the threshold where SageMath 5.6 invoked a $p$-adic algorithm to calculate determinants of integer matrices.  SageMath includes several different implementations of algorithms to compute integer determinants. The default algorithm depends on the input but can be explicitly overridden. Choosing a different algorithm made the computation work, so it was clear that the problem was with the $p$-adic algorithm.

How is the $p$-adic algorithm supposed to work?  Suppose we are given an $n \times n$ integer matrix $A$.  Checking whether the determinant of $A$ is 0 is relatively quick, so let us assume that $\det(A) \neq 0$.  The first step is to choose a random integer vector $v$ and a prime number $p$.  We may find a rational vector $x$ that solves the matrix equation $Ax = v$ using $p$-adic approximation.  That is, we solve $A \bar{x} = b \pmod{p^m}$ for some $m$; if we take $m$ to be sufficiently large, we may recover $x$ from $\bar{x}$ (see \cite{dixon}).  Now, one can prove, using Cramer's rule, that the least common multiple of the denominators of the entries of $x$
will be a divisor $d$ of $\det(A)$.  With high probability, $\det(A)/d$ will be not only an integer, but a tiny integer.  We can reconstruct $\det(A)$ completely by finding its value modulo a few different prime numbers, using the row-echelon form of $A$ in fields of prime order, and then applying the Chinese Remainder Theorem.  The Hadamard bound on the determinant of a matrix, which compares the determinant to the product of the Euclidean norms of its columns, bounds the size of the prime numbers we need to check: this guarantees that we can complete the reconstruction in a finite amount of time.  (Sage's implementation of the $p$-adic algorithm actually approximated the solution to $Ax=v$ using several primes $p_i$; details and time estimates may be found in \cite{iml}.)

Further testing of the $p$-adic determinant algorithm in SageMath showed that it worked for very large integer matrices, as well as small ones.  One of us, Jeroen Demeyer, a postdoc at Ghent University, discovered that in fact the computation only stalled for $51 \leq n \leq 63$.  He changed the title of the ticket to \emph{determinant() of integer matrices of size in [51,63] broken}, and began trying to figure out why SageMath treated these matrices differently.  The problem, he discovered, lay in the implementation of the last part of the $p$-adic algorithm, where SageMath tried to find $\det(A)$ ``modulo a few prime numbers''.  When a prime $p$ is large, SageMath computed determinants $\pmod{p}$ using the code for determinants over $\mathbb{Z}$.  A recent tweak to another part of Sage's matrix code had changed the definition of ``large $p$'' to be $p>2^{23}$ (that is, prime numbers greater than $8 388 593$).  When $n \leq 63$, Sage's integer matrix determinant function called code that asked for $\det(A)$ modulo primes greater than $8 388 593$.  This produced an infinite loop: computing the determinant of an integer matrix called for the determinant $\pmod p$, which called for the determinant of an integer matrix.

Once he found the bug's origin, Demeyer quickly wrote a patch to fix it.  The next step is peer review: before someone's code is incorporated in Sage, another developer must test it and sign off.  Review often entails multiple rounds of suggestions and fixes, but in this case Demeyer's fix worked well the first time, and the mathematical physicist Volker Braun (then a postdoc at the Dublin Institute for Advanced Studies) approved it for inclusion in Sage, just two days after the bug was reported. Less than a month later, SageMath 5.7, which included the fix for this bug, was released.

Part of Demeyer's patch was a test to ensure similar problems would not arise in future.  Each new release of SageMath goes through a series of tests to ensure that it works as advertised.  Demeyer's patch requires new versions of the SageMath code to compute the determinants of random integer matrices up to $80 \times 80$ successfully.  Furthermore, SageMath now checks that $\det(A^2) = \det(A)^2$ for a random square integer matrix $A$.  Automatic testing provides SageMath with resilience and makes it easier for multiple developers to contribute to different parts of the code base.  Currently, about 94\% of the functions in SageMath are automatically tested.

By default, new releases of SageMath use FLINT (\cite{flint}), the Fast Library for Number Theory, to compute the determinants of integer matrices.   For most integer matrices $A$ of size $24 \times 24$ and higher, FLINT computes the determinant $\det(A)$ by calculating $\det(A)$ modulo small primes.  The primes are chosen so that their product is greater than twice the Hadamard bound.  When the entries of $A$ are small in comparison to the size of the matrix, however, FLINT still uses the $p$-adic lifting algorithm.

\end{document}